\documentclass[11pt]{article}
\usepackage{sbc-template}
\usepackage{graphicx,url}
\usepackage[latin1]{inputenc}
\usepackage{graphics}
\usepackage{amssymb,amsthm,amsmath}
\usepackage[rm]{subfigure}    
\usepackage{threeparttable}   
\usepackage{algorithm}               
\usepackage{algorithmic}
\usepackage{multirow}
\usepackage[dvips]{color}
\usepackage{lineno}

\usepackage{tikz}
\definecolor{gold}{rgb}{0.85,.66,0}

\newtheorem{definition}{Definition}

\newtheorem{lemma}{Lemma}

\begin{document}

\sloppy

\title{ Energy and Spectral Efficiencies Trade-off with Filter Optimization in Multiple Access Interference-Aware}
\author{Álvaro R. C. Souza\inst{1}, Taufik Abrão\inst{1}, Lucas H. Sampaio\inst{2}, Paul Jean E. Jeszensky\inst{2}}
\address{Dept. of Electrical Engineering and Computer Science\\ State University of Londrina, 86051-970. Po.Box 6001, Brazil
\nextinstitute Dept. of Telecommunications and Control Engineering\\ Escola Politécnica of the University of São Paulo. Brazil
\vspace{-1mm}
\email{\{alvarorcsouza, lucas.dias.sampaio\}@gmail.com; \,\, abrao@ieee.org}
}

\maketitle

\begin{abstract}
This work analyzes the optimized deployment of two resources scarcely available in mobile multiple access systems, i.e., spectrum and energy, as well as the impact of filter optimization in the system performance. Taking in perspective the two conflicting metrics, throughput maximization and power consumption minimization, the distributed energy efficiency (EE) cost function is formulated. Furthermore, the best energy-spectral efficiencies (EE-SE) trade-off is achieved when each node allocates exactly the power necessary to attain the best SINR response, which guarantees the maximal EE. To demonstrate the validity of our analysis, two low-complexity energy-spectral efficient algorithms, based on distributed instantaneous SINR level are developed, and the impact of single and multiuser detection filters on the EE-SE trade-off is analyzed.
\end{abstract}
\vspace{-2mm}

\section{Introduction}
Resource allocation (RA) techniques, mainly power optimization, are becoming increasingly important in wireless system design, since battery technology evolution has not followed the explosive demand of mobile devices. The aim in RA is to maximize the sum of utilities of link rates for best-effort traffic. The usual approach consists in treat the problem jointly, i.e., optimize the joint power control and link scheduling, which has been extensively investigated in the literature and is known to be notoriously difficult to solve, even in a centralized manner. Hence, the methodology in \cite[Ch.4-6]{Boche_book06} consists in identify a class of utility functions for which the power control problem can be converted into an equivalent convex optimization problem. The convexity property is a key ingredient in the development of powerful and efficient power control algorithms.

One of the most interesting way of dealing with power allocation problem is the energy-efficiency (EE) approach \cite{Meshkati_07, Buzzi10, Miao10}, with aims to maximize the transmitted data per energy unit (measured in bits per Joule) and closely related to green communication techniques \cite{Han11}. As pointed out by \cite{Chen11}, one of the most important trade-offs on green wireless communications is energy efficiency \emph{versus} spectral efficiency trade-off (EE-SE); the goal consists in balancing these two important conflicting metrics.

Recently, game theory, which has its roots in the economy field, has been broadly applied to wireless communications for random access and power control optimization problems. This work proposes a power control procedure based on the optimized deployment of two main resources scarcely available at the multiple access mobile terminals (MT's), i.e., spectrum and energy. Importantly, from the analysis of two conflicting metrics, throughput maximization and power consumption minimization, the distributed energy efficiency cost function is formulated as a non-cooperative game. Indeed, the overall EE network depends on the behavior of each single user; so, the power control can be properly modeled as non-cooperative game \cite{Fudenberg91}.

This work also investigates the impact of multiuser detection schemes, motivated by the fact that the gap between optimal-EE and maximal-SE is reduced when the multiple-access interference (MAI) is increased. Since those detectors are capable to reduce substantially the MAI from other users, their deployment could result in more energy-efficient systems, meaning the same SINR can be achieved with less power consumption. In order to demonstrate the validity of the method, two low-complexity energy-spectral efficiency algorithms based on distributed instantaneous SINR level are developed.

\subsection{Related works}
The energy efficiency power allocation problem with filter optimization was developed for DS-CDMA systems in \cite{Meshkati05, Buzzi08, Buzzi10}, and demonstrate the impact of MuD strategies in the EE maximization. Additionally, \cite{Zappone11} introduced cooperative networks in that scenario. For multi-carrier systems, \cite{Miao09, Miao10} investigates the energy efficiency problem for OFDMA systems, and the EE-SE trade-off for OFDMA systems are investigated in \cite{Miao11}. Considering multi-carrier CDMA, \cite{Meshkati06} investigates the existence of Nash Equilibrium for the EE optimization problem, and even assuming some simplifications the system can achieve multiple or even none Nash Equilibrium. To reduce the allocated power by non-cooperative games, \cite{Saraydar02} develops a linear pricing factor, and demonstrate that this approach is Pareto-dominant over non-pricing solutions.

This paper proposes the EE-SE trade-off analysis in DS-CDMA systems, as developed by \cite{Miao11} for OFDMA systems. Based on the conclusions of this analysis, we develop two algorithms to improve the EE of the system by removing EE non-optimal users, putting these users in outage. So far, this approach was not presented in other works. Hence, in previous literature's works, when a user is not able to achieve the optimal EE, it uses the maximum power to achieve the maximum EE. 
Since CDMA systems are limited by interference, removing those users it enables to reduce the generated interference and, possibly, increases the energy efficiency.

\section{Network System Model} \label{sec:sys_model}
For analysis simplicity, initially we have assumed a single rate uplink direct sequence code division multiple access (DS/CDMA) network. However, the extension for multi-cell multi-carrier multiple access systems is straightforward. The received signal in the base station can be described as:
\begin{equation}
    \label{eq:sinal_bs}
    \mathbf{y} = \sum_{k=1}^{K}\sqrt{p_k}h_{k}b_{k}\mathbf{s}_k + \pmb{\eta}
\end{equation}
\noindent where $h_{k}$ is the complex channel gain between the $k$th user and the base station, constant during the chip period\footnote{Mobile channel is assumed to be slow and non-selective in frequency.}, $\mathbf{s}_k$ is the $k$th user spreading code with length $N$, representing the processing gain; the modulated symbol is given by $b_k$, and $\pmb{\eta}$ is the thermal noise, assumed to be AWGN, zero-mean and covariance matrix given by $\sigma^2\mathbf{I}_N$.

The uplink $1 \times K$ channel gain vector, considering path loss, shadowing and fading effects, between users and the base station, is given by:
\begin{eqnarray}
{\bf h}_{\rm upl} = \left[ h_1, \cdots, h_K \right]
\end{eqnarray}
\noindent which could be assumed static or even dynamically changing over the optimization window ($N$ time slots). The signal-to-interference-plus-noise ratio (SINR) is defined by the received signal power to the sum of interfered power plus background noise, measured after demodulation. In DS/CDMA this ratio depends on the detection strategy. Considering the adoption of linear receivers, the SINR can be expressed, generically, as:
\begin{equation}\label{eq:sinr_gen}
\gamma_k = \frac{p_k |h_{k}|^2 ({\bf d}_k^T {\bf s}_k)^2}{\sum\limits_{j=1 \atop j \neq i}^{K} p_j |h_{j}|^2 ({\bf d}_k^T {\bf s}_j)^2 + \sigma_k^2 ({\bf d}_k^T {\bf d}_k)} = \frac{p_k |h_{k}|^2}{I_k+\sigma_k^2 ({\bf d}_k^T {\bf d}_k)}
\end{equation}
\vspace{-3mm}

\noindent where the channel gain $h_i=|h_{i}|e^{\angle h_{i}}$, $I_k = \sum_{j \neq i}^{K} p_j |h_{j}|^2 ({\bf d}_k^T {\bf s}_j)^2$ represents the MAI power level, ${\bf s}_k = \frac{1}{\sqrt{N}}[c_1, c_2, \cdots, c_N], c_i = \mathcal{U}\{-1,1\}$ is the $k$th user pseudo-noise (PN) spreading code, with $({\bf s}_k^T {\bf s}_k) = 1$, and ${\bf d}_k$ is the $N-$dimensional vector representing the receive filter for the $k$th user;  $(\cdot)^T$	denotes transpose operator.

\subsection{Matched Filter (SuD)}\label{sec:sud_filter}
The simplest filter that can be used is the matched filter (MF), a single-user detection (SuD) strategy. For this receiver, the filter vector ${\bf d}_k$ is defined as the $k$th user spreading code, and the interference power is considered as a background noise, which limits the system performance, since CDMA systems are limited by the interference level. Hence, the SINR expression in eq. \eqref{eq:sinr_gen} can be re-written considering MF by:
\begin{equation}
    \label{eq:sinr_mf}
    \gamma_k^{\rm MF} = \frac{p_k |h_{k}|^2 ({\bf s}_k^T {\bf s}_k)^2}{\sum\limits_{j=1 \atop j \neq i}^{K} p_j |h_{j}|^2 ({\bf s}_k^T {\bf s}_j)^2 + \sigma_k^2 ({\bf s}_k^T {\bf s}_k)} = \frac{p_k |h_{k}|^2}{\sum\limits_{j=1 \atop j \neq i}^{K} p_j |h_{j}|^2 ({\bf s}_k^T {\bf s}_j)^2 + \sigma_k^2} = \frac{p_k |h_{k}|^2}{I_k^{\rm MF} + \sigma_k^2}
\end{equation}

Besides the simplicity, consider the MAI as noise power implies that the system cannot mitigate the interference, and when the system loading becomes higher, the necessary power to keep the SINR level for a specific user becomes higher; as a consequence, this impacts on the whole multiple access system, increasing the overall power level consumption. In order to avoid this behavior, improve the system performance and simultaneously reduce the MAI effect, Verdu developed the idea of multi-user detection (MuD) \cite{Verdu84}. In order to detect/decode the signal of the interest user, multiuser receivers deploy MAI information (come from interfered signal users), reducing the necessary power to achieve the same SINR level. The best MuD strategy namely optimum MuD receiver (OMuD) is that this receiver results in an exponential complexity, which reduces the applicability of the optimum receiver.

\subsection{Sub-Optimum Linear Multiuser Filters (LMuD)}\label{sec:mud_filter}
One of the possible ways to reduce the OMuD complexity, obtaining near-optimal performance, consists in apply linear multiuser filters \cite{Lupas89}, such as decorrelator (DE), zero-forcing (ZF) and the minimum mean square error (MMSE). Linear multiuser filters applies a linear transformation in the MF soft estimation, decoupling the MAI. Besides, LMuDs are useful in power optimization algorithms because the resulting SINR is deterministic, unlike heuristic-based methods, and that is important to find the minimum power to achieve the target SINR.

Among the LMuD filters, the most efficient is the MMSE, because this technique takes into account the amplitude and background noise from the interfered users, which results in a most efficient interference reduction without a large increase in the background noise, as occur for zero-forcing and decorrelator as well \cite{Moshavi96}. Besides the efficiency, the MMSE filter requires users' amplitude matrix, and a distributed implementation becomes complicated -- mainly whit the proposed iterative power control algorithm, as described in Section \ref{sec:pseudocodes}, since this matrix will be updated at each iteration.

In order to avoid the aforementioned problems, decorrelator detector has been chosen, which presents slight inferior performance regarding the MMSE detector, but it depends only on the spreading codes (${\bf s}_k$) and the correlation matrix ($\mathbf{R}$) \cite{Meshkati05}. Both parameters are constant during the power control algorithm execution, which implies in just one transmission at the beginning of algorithm iterations. The decorrelator filter (after MF) is given by:
\begin{equation}
    \label{eq:filter_dec}
    {\bf d}_{\rm DEC} = [{\bf d}_1\, {\bf d}_2\, \cdots\, {\bf d}_k\, \cdots\, {\bf d}_K] = {\bf S} ({\bf S}^T{\bf S})^{-1} = {\bf S}\mathbf{R}^{-1}
\end{equation}
\noindent Hence, the achieved SINR is given by:
\begin{equation}
    \label{eq:sinr_dec}
    \gamma_k^{\rm DEC} = \frac{p_k |h_k|^2}{\sigma^2 {\bf d}_k^T {\bf d}_k} = \frac{p_k |h_k|^2}{I_k^{\rm DEC}}
\end{equation}

The investigation of the tradeoff between the resource reduction achieved with the MMSE detector (instead of decorrelator), and the needed power and communication overhead to transmit the amplitude information is out of the scope of this work and could be addressed in a future work.

\subsection{QoS Requirements}
In order to guarantee the quality of service (QoS),  a minimum data rate $R_{k,\min}$ for each user must be provided for the system network service, being an important requirement to be warranted. So, in general, data rate for the $k$th user is assumed to be a function of SINR $\gamma_k$. To do that, we use a modified version of Shannon capacity equation, given by:
\begin{equation}
    \label{eq:shannon2}
    r_k = \mathcal{C}_k^{\rm gap} = \mathtt{w} \log_2 (1+ \theta_k \cdot \gamma_k), \qquad \forall k \qquad {\small \left[\rm bit / s\right]}
\end{equation}
where $\gamma_k^{\rm filter}$ is given herein by the \eqref{eq:sinr_mf} or \eqref{eq:sinr_dec} and $\theta_k$ is a gap introduced to describe the limitations and imperfections in real communication systems, such as modulation effects, and so on, to approximate the real data rate \cite{TSE_05}, given by:
\begin{equation*}
    \theta_k = - \dfrac{1.5}{\log(5\, \textsc{ber}_k)}, \qquad \text{with}\,\,\ \theta_k\in [0;1[
\end{equation*}
\noindent and $\textsc{ber}_k$ is the maximum tolerable bit error rate by the $k$th user \cite{Gold_97}.
The spectral efficiency (SE) is obtained from (\ref{eq:shannon2}):
\begin{equation}
    \label{eq:SE}
    \eta_k = \log_2 (1+ \theta_k \cdot \gamma_k), \qquad \forall k \qquad {\small \left[\rm \frac{bit}{s \cdot Hz} \right]}
\end{equation}

From \eqref{eq:shannon2}, the minimum data rate for the $k$th link, $R_{k,\min}$, which is able to guarantee the QoS, considering maximum tolerable BER for that service, can be easily mapped into the minimum SINR:
\begin{eqnarray}
\label{eq:newshannon2}
\gamma_{k,\min} &=& \frac{2^{\frac{R_{k,\min}}{r_c}} - 1}{\theta_k} \qquad \forall k = 1,\ldots,K
\end{eqnarray}

\vspace{-6mm}
\section{Problem Formulation}\label{sec:problem}
In a multiple access interference limited communication system, the $k$th user selfish (non cooperative approach) allocates his own transmit power $p_k$ and receive filter strategy (single- or multi-user detection strategy, cancelation, MAI mitigation, zero-forcing and so forth) in order to maximize his own energy efficiency function, expressed by \cite{Goodman00}:
\begin{equation}
    \label{eq:util_basic}
    \xi_k  = r_k \frac{L}{M}\frac{f(\gamma_k)}{p_k + p_c} \qquad \left[\frac{\rm bit}{\rm Joule}\right], \qquad \forall k=1,\ldots, K
\end{equation}
\noindent where $M$ is the number of bits in each transmitted data packet; $L$ is the number of information bits contained in each data packet, $p_k$ is the transmission power, $p_c$ is the circuit power consumption, and $f(\gamma_k)$ is the efficiency function, which approximates the probability of error-free packet reception. When no coding technique is used, it can be approximated by
\begin{equation}
    \label{eq:eff_func}
	f(\gamma_k)=(1-e^{-\gamma_k})^M
\end{equation}
\noindent This approximation is widely accepted for BPSK and QPSK modulation.

It is worth noting that both transmission power and circuit power consumptions are very important factors for energy-efficient communications. While $p_k$ is used for reliable data transmission, circuit power represents average energy consumption of electronics devices and circuitry \cite{Miao10}. Besides, the SINR for user $k$, $\gamma_k$, assumes different definition depending on system type, multiple access detection strategy (SuD or MuD, as described in Section \ref{sec:sud_filter} and \ref{sec:mud_filter}, respectively), spreading sequence type an so forth.

Note that $\xi_k$ is measured in $\left[\frac{\rm bit}{\rm Joule}\right]$, which represents the number of successful bit transmissions that can be made for each energy-unit drained from the battery and effectively used for transmission.

In a more general context, we can define the concept of global energy efficiency function as the ratio of the total achievable capacity over the total power transmission consumption:
\begin{equation}
    \label{eq:EE_se}
    \bar{\xi}  =  \frac{\sum_{k=1}^{K}\ell_k r_k  f(\gamma_k)}{P_{\rm Tot}}\qquad \left[\frac{\rm bit}{\rm Joule}\right]
\end{equation}
\noindent where $P_{\rm Tot}= \sum_{k=1}^K (p_k + p_c)$,  and $\ell_k=\left(\frac{L}{M}\right)_k$

\subsection{Distributed Non-cooperative EE Power Optimization Game}
The network energy efficiency depends on the behaviors of all users; so, the power control can be properly modeled as non-cooperative game \cite{Fudenberg91}. In the context of non-cooperative power control game
\begin{equation}
    \label{eq:game}
	\mathcal{G} = \left[\mathcal{K}, \, \left\{\mathcal{A}_k\right\}, \, \left\{u_k\right\}\right]
\end{equation}
\noindent where $\mathcal{K}= \left\{1,2,\ldots, K \right\}$, and $\left\{\mathcal{A}_k\right\}=[0, P_{\max}]$ is the strategy set for the $k$th user, with $P_{\max}$ being the maximum allowed power for transmission; the utility functions $\left\{u_k\right\}$ is performed by one of the energy efficiency functions $\left\{\xi_k\right\}$, such as the basic EE function, eq. (\ref{eq:util_basic}).


Consider the power allocation for the $k$th user, $p_k$ and denote the respective power vector of other users (interfering users):
\begin{equation}
{\bf p}_{-k} = [p_1, p_2,\ldots, p_{k-1}, p_{k+1}, \ldots,p_{K}]
\end{equation}
Hence, given the power allocation of all interfered users, $p_{-k}$, the best response of the power allocation for the $k$th user can be  expressed as:
\begin{equation}
\label{eq:pbest}
p_{k}^{\rm best} = f_k({\bf p}_{-k}) = \arg \max_{p_k} \,\,\, u_k(p_k,{\bf p}_{-k})
\end{equation}
where $u_k$ is given by (\ref{eq:util_basic}), and $f_k({\bf p}_{-k})$ is called the $k$th best response function.

Finally, the problem for distributed energy-efficiency with power constraint under non-cooperative game perspective can be posed as:
\begin{eqnarray}\label{eq:ee_prob}
\arg\, \max_{p_k} \, \xi_k &=& \arg\, \max_{p_k} \,\,\,\, \ell_k r_k \frac{f(\gamma_k)}{p_k+p_c}\\
	                     s.t. && 0 < p_k\leq P_{\max}\nonumber
\end{eqnarray}
\noindent which solution consists in adopting the best-response strategy for user $k$. Indeed, the  best-response strategy consists in obtain the best user utility function (EE) individually for each user, as posed by (\ref{eq:pbest}). Hence, applying the derivative on eq. (\ref{eq:ee_prob}), regarding $p_k$, we have:
$$
\frac{\partial \xi_k}{\partial p_k} = 0
$$
which, under certain conditions, represents the best response power allocation strategy for each user, given the interfered power vector ${\bf p}_{-k}$.

\subsection{Best SINR Response for SuD and LMuD Filters}
Given the context of medium or high SINR, the power allocated for $k-$th can be approximated by
\begin{equation}
	p_k \approx \gamma_k \frac{I_k}{|h_k|^2} = \gamma_k \widetilde{I}_k^{\rm MF}\ ({\rm MF}), \qquad p_k = \gamma_k \frac{I_k}{|h_k|^2} = \gamma_k \widetilde{I}_k^{\rm DEC}\ ({\rm DEC}).
\end{equation}

The first derivative of EE function (\ref{eq:ee_prob}) regarding to $p_k$ is equivalent to take the derivative of EE function regarding to $\gamma_k$:
\begin{equation}\label{eq:EE_deriv}
\frac{\partial \xi_k}{\partial \gamma_k} = \frac{\partial}{\partial \gamma_k}
 \left\{\ell_k \frac{(1-e^{-\gamma_k})^M \log(1+\theta_k\gamma_k)}{\gamma_k \widetilde{I_k} + p_c} \right\}
\end{equation}

Hence, the optimal SINR for the $k$th user, $\gamma_k^*$, in terms of EE-SE tradeoff is obtained finding the solution of $\frac{\partial \tilde{\xi}_k}{\partial \gamma_k} =0$ (maximization point), admitting fixed the normalized multiple access interference $\widetilde{I}_k$. This condition is equivalent to solve function (\ref{eq:ee_se_deriv}) regarding $\gamma_k$.
\begin{equation}\label{eq:ee_se_deriv}
Me^{-\gamma_k} \, \log_2(1+\theta_k\gamma_k) + \frac{\theta_k (1-e^{-\gamma_k}) }{(1+\theta_k\gamma_k)\ln2}  =   \frac{\widetilde{I_k} \log_2(1+\theta_k\gamma_k) (1-e^{-\gamma_k})}{(\gamma_k\widetilde{I_k}+p_c)}
\end{equation}

In order to guarantee that eq. \eqref{eq:ee_se_deriv} has only one maximizer, we introduce the concept of quasiconcavity, defined as \cite{Miao11}:
\vspace{3mm}

\begin{definition}[Quasiconcavity]
    A function $z$, that maps a convex set of n-dimensional vectors $\mathcal{D}$ into a real number is quasiconcave if for any $\mathbf{x}_1, \mathbf{x}_2 \in \mathcal{D}, \mathbf{x}_1 \neq \mathbf{x}_2$
    \begin{equation}
            \label{eq:quasiconcavity}
            z(\lambda \mathbf{x}_1 + (1 - \lambda)\mathbf{x}_2) \geq \min\left\{ z(\mathbf{x}_1), z(\mathbf{x}_2) \right\}, \qquad \text{where},\,\, \lambda \in (0,1).
    \end{equation}
\end{definition}

\begin{lemma}[Quasiconcavity of $u_k$]
    The utility function $u_k(p_k, \mathbf{p}_{-k})$ is quasiconcave in $p_k$
    \label{lemma:quasiconcavity}
\end{lemma}
This result is very important in the proof of existence and uniqueness of the system equilibrium. However, due to space limitation all proofs are not developed herein.

\vspace{-4mm}
\section{Increasing Interference Effect and Nash Equilibrium on EE-SE Trade-off} \label{sec:interf_tradeoff}
\begin{figure}[ht]
\begin{minipage}[b]{0.44\linewidth}
\centering
    \includegraphics[width=1\textwidth]{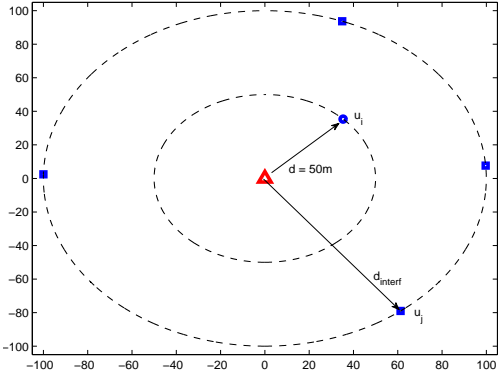}
    \caption{a) Cell geometry with increasing interference level, $I\varpropto d_{\rm interf}^{-1}$. b) Adopted System Parameters for EE-SE Trade-off analysis}
\label{fig:cell_geo}
\vspace{-35mm}
\end{minipage}
\vspace{-5mm}
\begin{minipage}[b]{0.57\linewidth}
	\centering
\footnotesize
\begin{tabular}{ll}
\hline
\textbf{Parameters} & \textbf{Adopted Values}\\
\hline
\hline
\multicolumn{2}{c}{\textit{DS/CDMA Optimal Power Allocation}} \\
\hline
Noise Power   & $P_n=-90$ [dBm]\\
Processing Gain  & $N = 15$\\
Max. power per user & $P_{\max}=10$ [dBm]\\
\# mobile terminals   & $K \in \{3; \, 9\}$\\
Interest user distance   & $d = 50$ [m]\\
Interfering users distance  & $d_{int} = [80, 100, 200]$ [m]\\
Packet size & $M = 80$ [bits] \\
Data bits & $L = 50$ [bits] \\
Maximum BER & $\mathrm{BER}_{k} = 10^{-3}$\\
Circuit Power & $p_c = 7$ [dBm] \\
Bandwidth & $\mathtt{w} = 10^6$ [Hz] \\
\hline
\multicolumn{2}{c}{\textit{Channel Gain}} \\
\hline
Path loss  & $\varpropto d^{-2}$ \\
Fading coefficients & Rayleigh distribution\\
                    & mean over 5000 samples \\
\hline
\multicolumn{2}{c}{\textit{Verhulst PCA}} \\
\hline
Convergence factor  & $\alpha = 0.5$ \\
\# iterations  & $N_{it} = 500$\\
\hline
\end{tabular}
\end{minipage}
\end{figure}

\begin{figure}[!htb]
\centering
\subfigure{
   a)\includegraphics[width=0.455\textwidth]{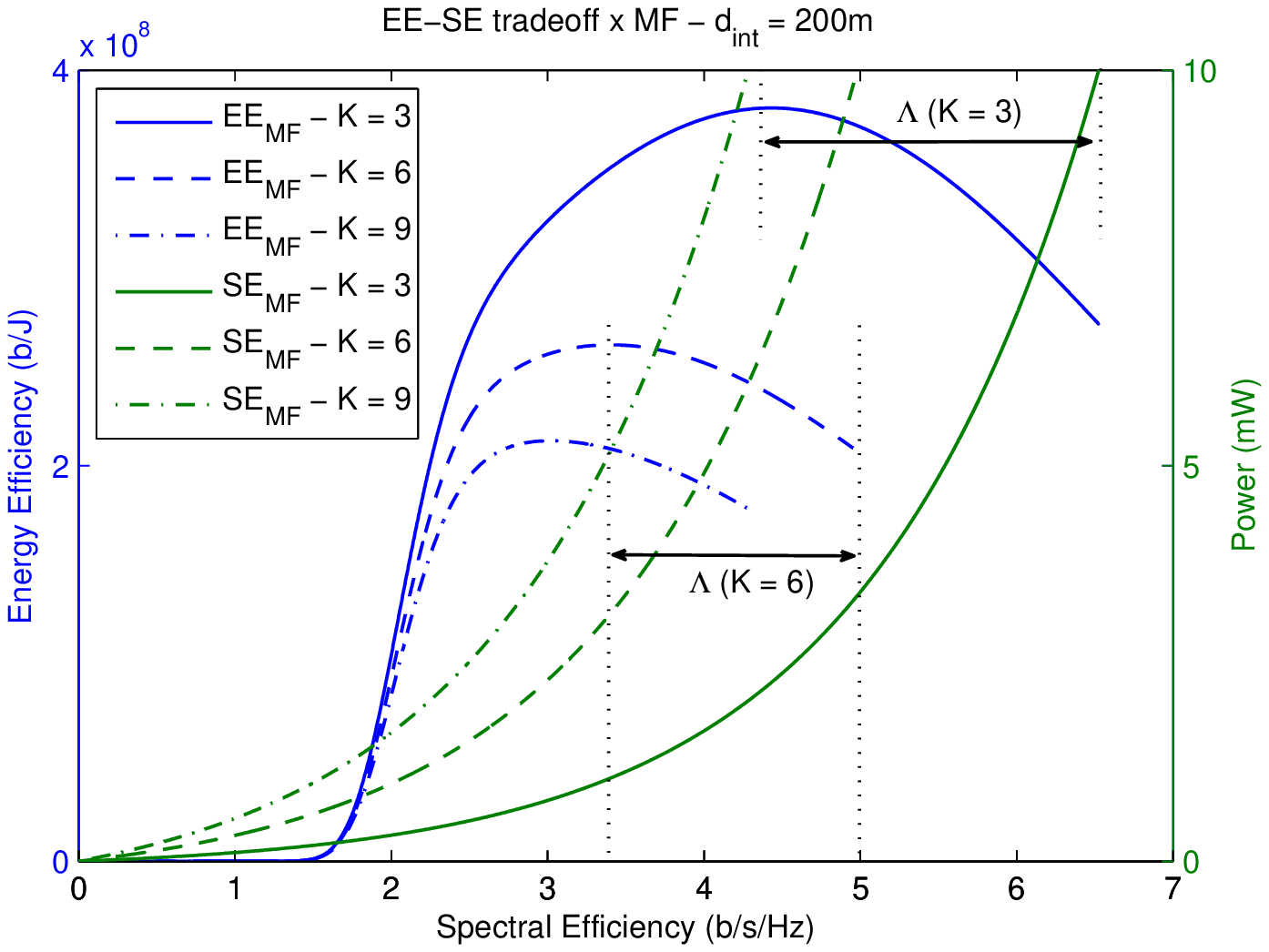}
     \includegraphics[width=0.455\textwidth]{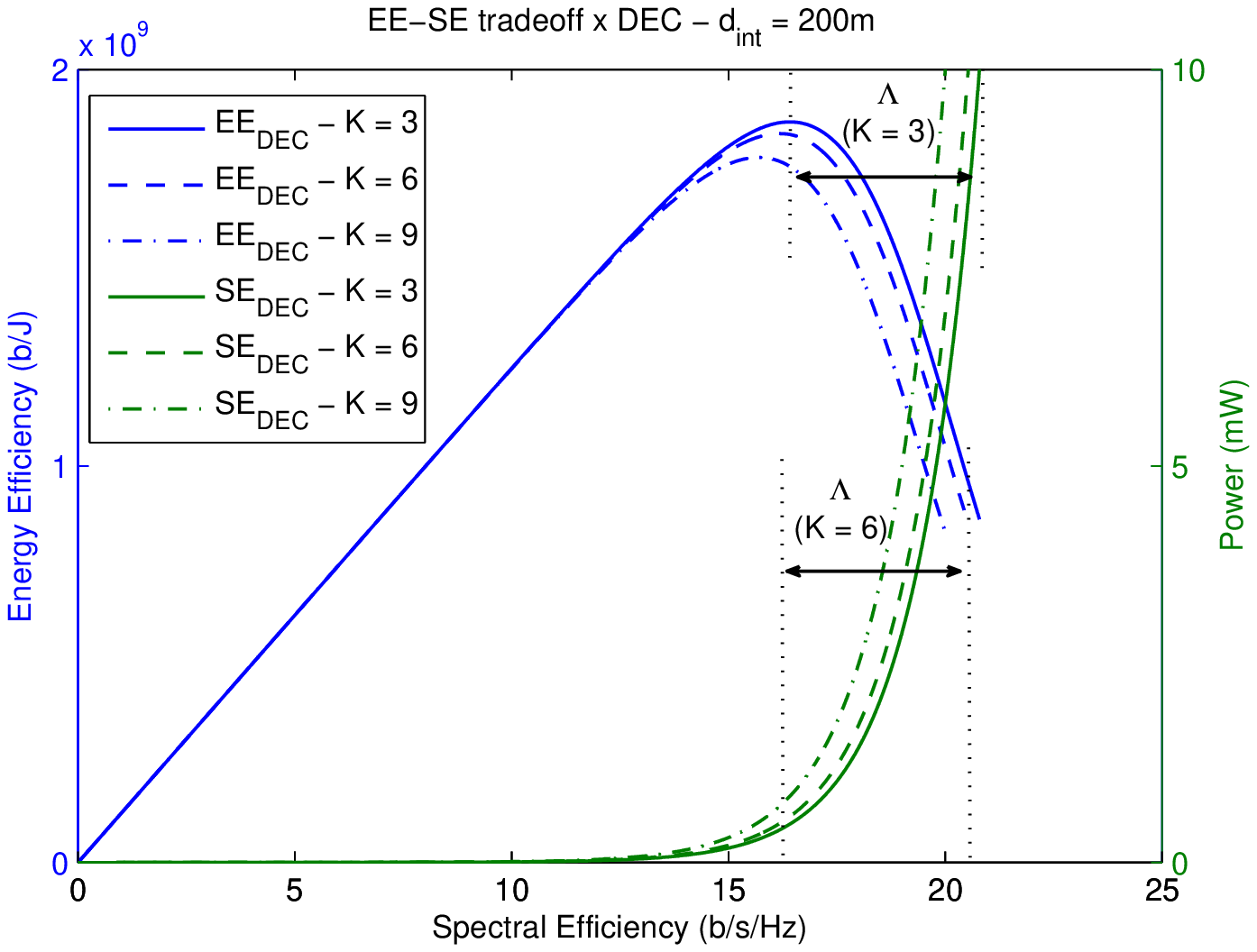}
 }\\
 \subfigure{
   b) \includegraphics[width=0.455\textwidth]{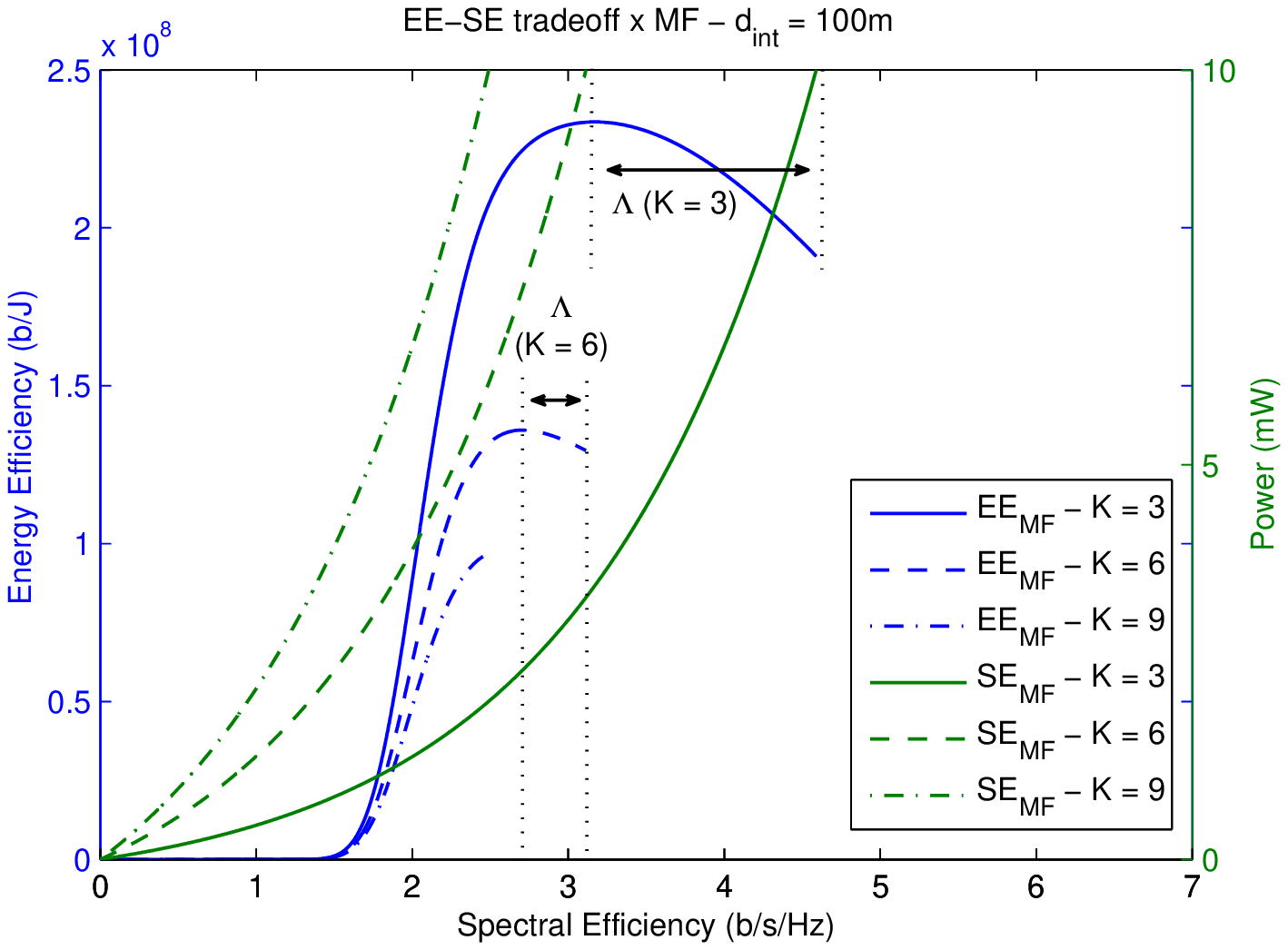}
   \includegraphics[width=0.455\textwidth]{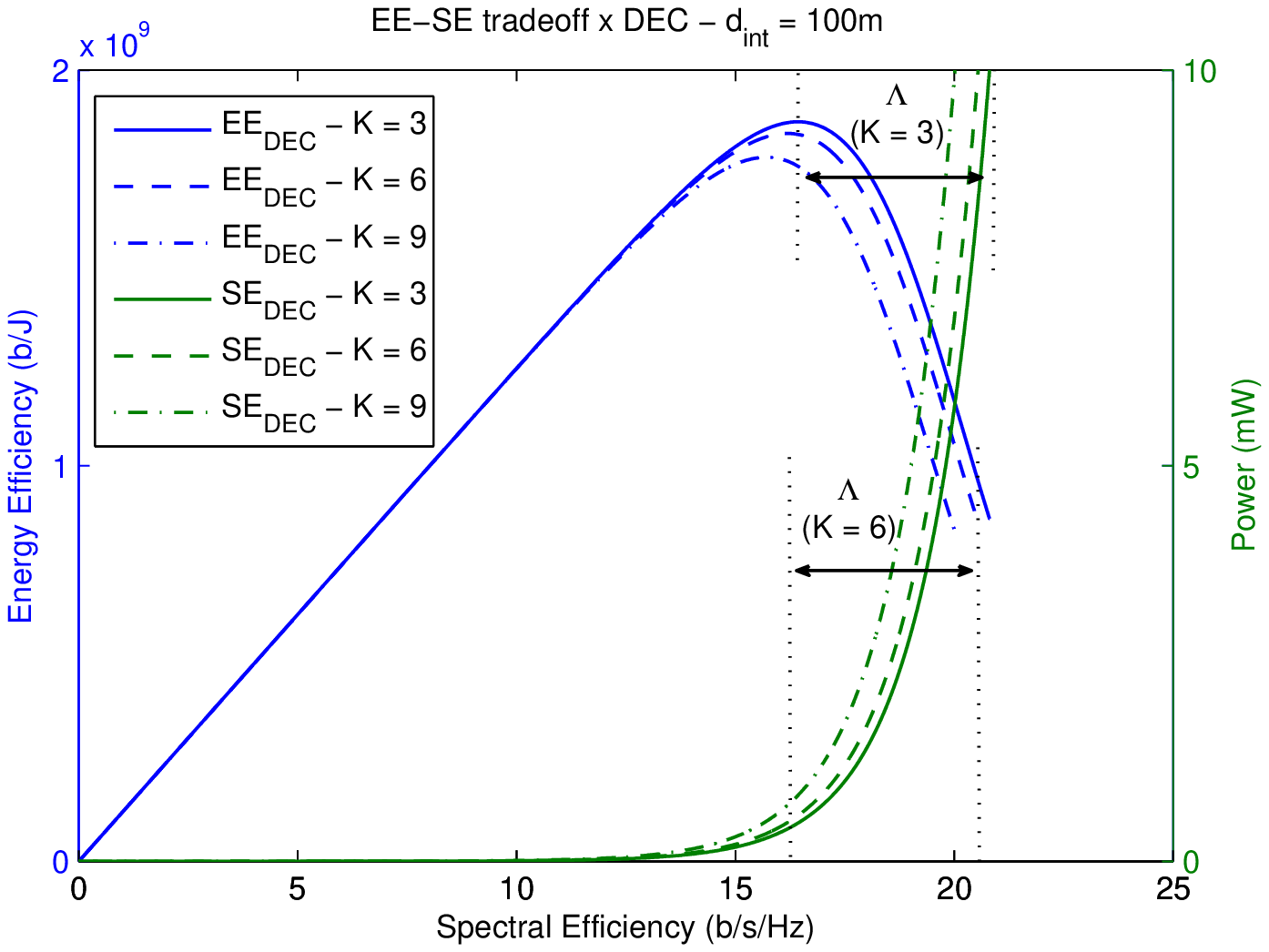}
 }\\
 \subfigure{
   c) \includegraphics[width=0.455\textwidth]{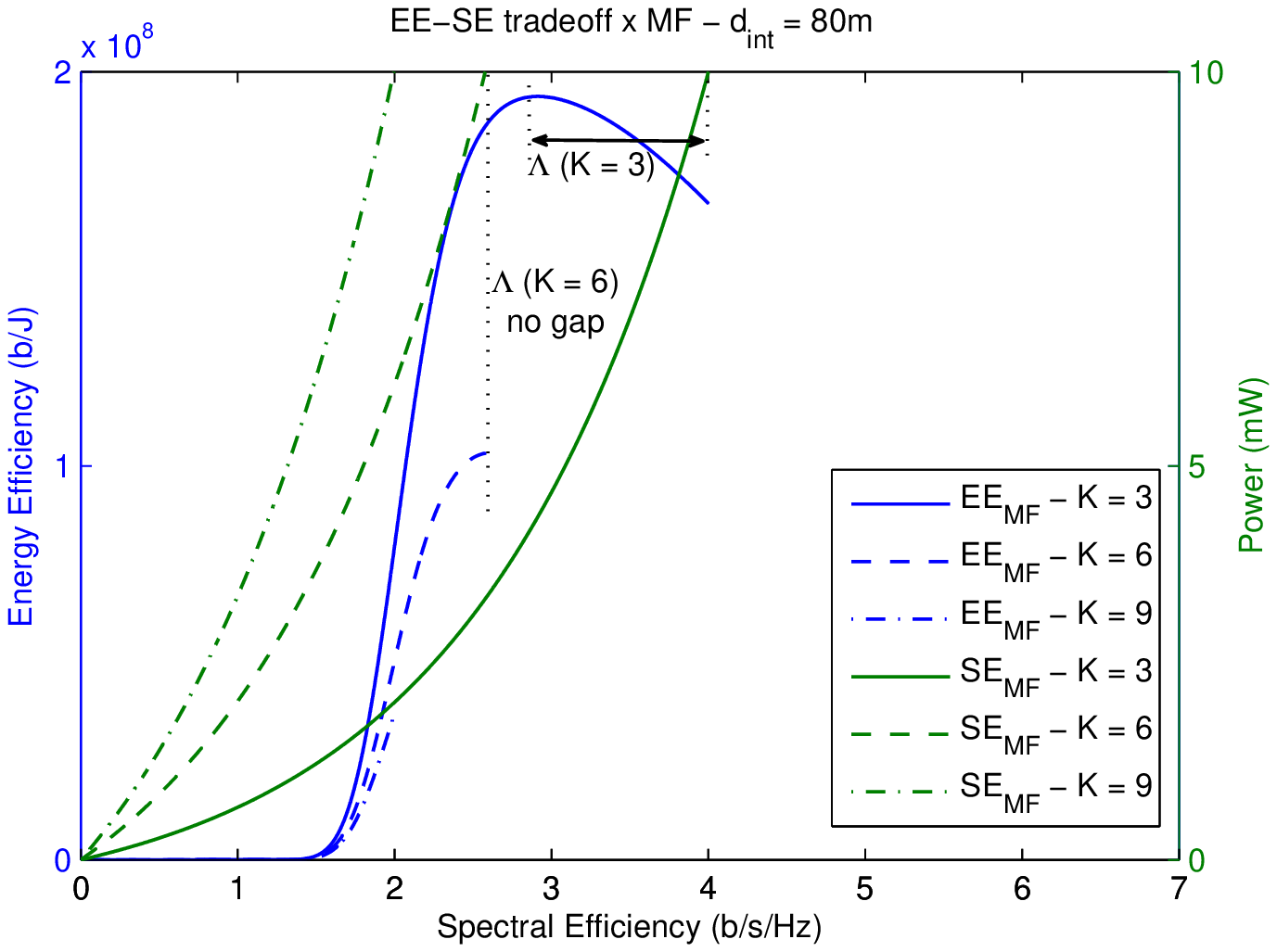}
   \includegraphics[width=0.455\textwidth]{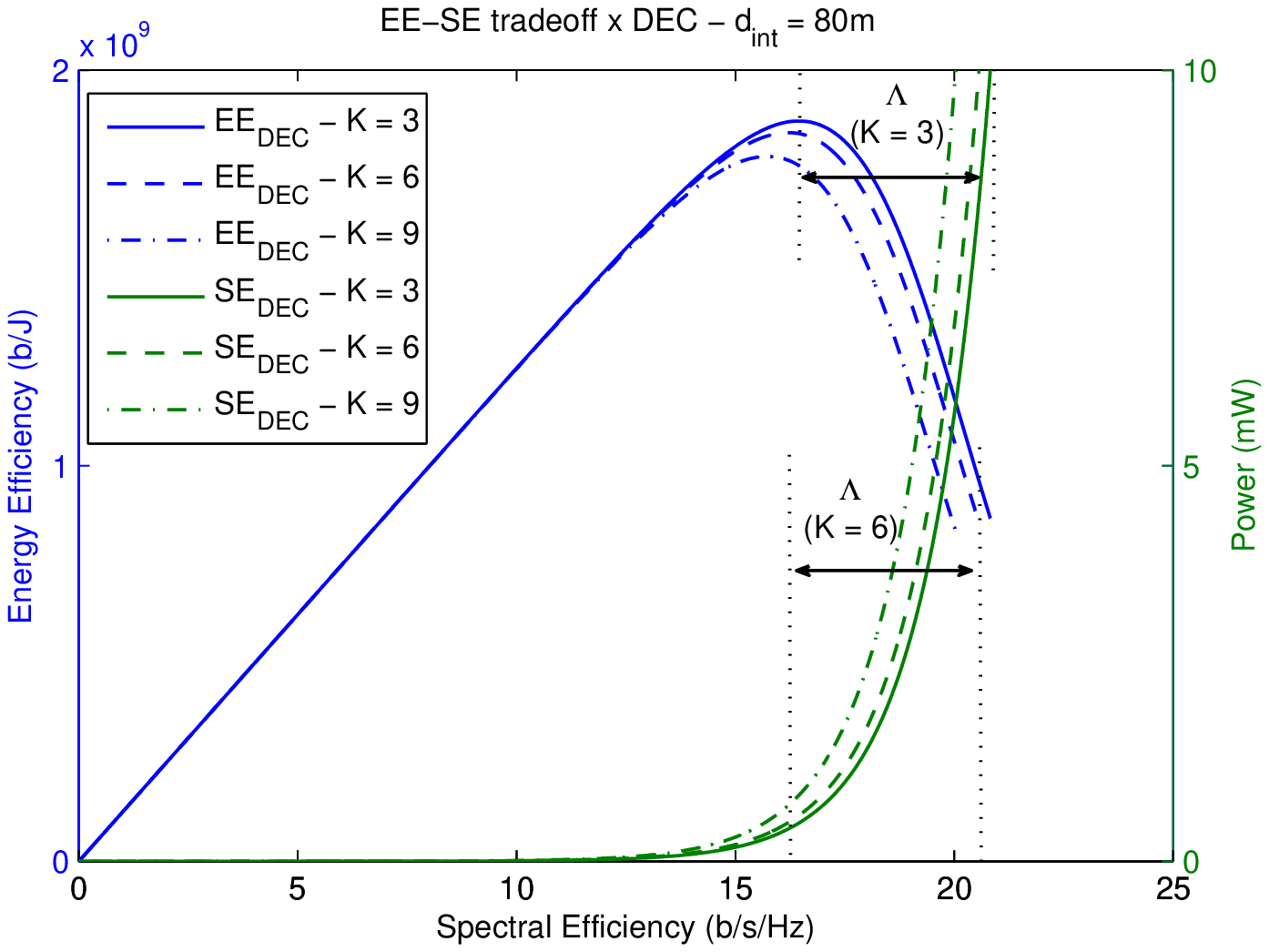}
 }
\vspace{-4mm}
\caption{EE-SE trade-offs considering different interfering scenarios and filters (left hand-side MF, right hand-side DEC). a) $d_{\rm interf}=200$m, $\beta_k = 0.25$; b) $d_{\rm interf}=100$m, $\beta_k = 0.50$; c) $d_{\rm interf}=80$m, $\beta_k = 0.63$.}
\vspace{-3mm}
\label{fig:ee-se_tradeoff}
\end{figure}

In this section we present a trade-off between non-cooperative energy-efficient and spectral-efficient power control schemes. This trade-off is determined by the multiple access interference level, which is responsible by the \emph{gap} among the maximal EE and the optimum SE (only attainable with infinity power allocation). In realistic interference-aware systems the increasing number of active users brings an increasing on system capacity; therefore, the SE of the system increases accordingly. The gap among the max-EE and the opt-SE, $\Lambda$, can be reduced when the interference level increases. In order to quantify this effect, let us define the coupling network parameter:
$$
\beta_k = \frac{\langle|h_k|^2\rangle}{\langle|h_j|^2\rangle},\qquad k: \text{interest};  \qquad j\neq k: \text{interfered users}
$$
where $\langle\cdot\rangle$ is the operator temporal average. Furthermore, defining the cell geometry and the placement of the $j=1,\ldots, K-1$ interfering users, as well as the $i$th interest user in the way of Fig. \ref{fig:cell_geo}.a, the max-EE and the opt-SE behavior are obtained in terms of $d_{\rm interf}^{-1}$,  Fig. \ref{fig:ee-se_tradeoff}. It is clear the gap reduction between the max-EE and the asymptotic-SE when the interference level increases, considering conventional detector. For decorrelator, the main components are the distance (since increases path loss) and spreading-code correlation, since bigger correlation implies in bigger noise amplification. Fig. \ref{fig:cell_geo}.b shows the system parameters used in this simulation scenario.

\vspace{-3mm}
\vspace{2mm}

When circuit power consumption is much smaller than the transmitted power ($p_c << p_k$), an interesting result emerges: the optimum SINR obtained from the EE optimization problem in \eqref{eq:ee_se_deriv} is the same for any MAI level, while the asymptotic SINR necessary to the SE maximization still remains related to the interference power level, $\widetilde{I}_k$.

Hence, under this hypothesis, the best SINR for max-EE criterium depends only on the system parameters, such as  maximal tolerable BER (QoS), modulation level, coding and packet coding size. It is worth to note that when the MAI increases, the transmitted power becomes higher and, indeed, the condition $p_c << p_k$ holds, as one can see from the left side plots of Fig. \ref{fig:ee-se_tradeoff}.a) to c), i.e., when distance $d_{\rm interf}$ is reduced.

In order to corroborate those conclusions and to determine the impact of the MAI in the energy efficiency problem, this work analyzes the impact of linear multiuser filter deployment, specially the use of decorrelator multiuser detector. Hence, the achieved SINR with decorrelator filter, eq. \eqref{eq:sinr_dec}, does not depend on the interference level, but the equilibrium point is changed by the position of the mobile terminals (i.e. the path-loss, which impacts on the channel gain, $h_k$), the instantaneous fading and the active users' spreading codes correlation (as one can be seen on the right side plots of Fig. \ref{fig:ee-se_tradeoff}). This way, the expectation is that the EE is almost the same for any system loadings when multiuser filter is deployed at receiver side.

\section{Proposed EE-SE Algorithms} \label{sec:pseudocodes}
The proposed algorithm to implement the optimal EE-SE tradeoff solution is described in Algorithm \ref{pcode:EE-SE-Verhulst} and is based on Verhulst power control algorithm \cite{TJG_10}. On the other hand, in order to avoid users' outage, in which users are not able to achieve the optimal SINR in terms of EE (due to $P_{\max}$ constraint), but are able to maintain the minimum data rate, $R_{k,\rm min}$, an alternative algorithm is proposed in Algorithm \ref{pcode:EE-SE-Verhulst_rate}.

Algorithm \ref{pcode:EE-SE-Verhulst} and \ref{pcode:EE-SE-Verhulst_rate} are closely related. It's easy to see that when the necessary SINR to achieve the minimum rate criterion $r_{k,\min}$ (and the associated $\gamma_{k,\min}$) is greater or equal than the optimum SINR ($\gamma_k^*$), Algorithm \ref{pcode:EE-SE-Verhulst_rate} reduces to Algorithm \ref{pcode:EE-SE-Verhulst}, since the second condition to be inserted in $K_{\rm out}$ is always true for the non-optimum users.

After defining the algorithms, we need to investigate the existence and uniqueness of the achieved equilibriums. Given that the equilibrium is defined by $\mathbf{p}^* = (p_1^*, p_2^*, \cdots, p_k^*)$, the Nash equilibrium can be defined as:
\vspace{2mm}

\begin{definition}[Nash Equilibrium]
    An equilibrium is said to be a Nash equilibrium if and only if any user cannot unilaterally improve their response by changing the optimum value. In the context of the energy-efficiency problem, this affirmation is equivalent to the fact that any user cannot improve their utility value by changing the optimum power for any other value:
\begin{equation}
u_k(p_k^*, \mathbf{p}_{-k}^*) \geq u_k(p_k, \mathbf{p}_{-k}^*), \quad  \forall k
        \label{eq:nash_eq}
\end{equation}
\label{def:nash_eq}
\end{definition}
\vspace{-6mm}

The uniqueness of the Nash equilibrium for both non-cooperative games is summarized in Lemma \ref{lemma:uniqueness}.

\begin{lemma}\label{lemma:uniqueness}
When the equilibrium $\bf{p}^*$ is achieved without removing any user, this Nash equilibrium is unique. When is needed to remove any user, multiple equilibriums will exist, depending on the adopted criterion. For our adopted criterion, the equilibrium is also unique.
\end{lemma}
\vspace{-2mm}

\section{Numerical Results} \label{sec:num_res}
System parameters are indicated previously in Fig. \ref{fig:cell_geo}.b with some changes. Analysis in this section assumes a ring geometry, with internal radius  $r_{\rm int}=50$m and external radius $r_{\rm ext}=200$m, with $K$ mobile users uniformly distributed on this ring area with radius $\sim \mathcal{U}[r_{\text{int}},\, r_{\text{ext}}]$, and the base station in the center of the ring. The processing gain was assumed $N = 63$; number of mobile terminals  $K \in \{2; \, 15\}$. For simplicity, identical parameters of QoS were adopted for all users, i.e., maximal tolerable $\mathrm{BER}_{k}=10^{-3}$, and minimum data rate $R_{\min} = 500 $ [kbps]. Fading is modeled as flat Rayleigh distribution (module), simulated by a complex Gaussian random process, with zero mean and variance given by $d^2_j$. In order to analyze the average network behavior, numerical results were taken as the average over 2000 network realizations, including random (uniform) users distribution, pseudo-noise spreading codes and channel gains; but all algorithms and filters have deployed the same statistical parameters. Furthermore, it was assumed that the mobile transmitter has perfect channel state information (CSI) available, but the measurement of other mobile users' CSI can only be carried out by the base station through quantized transmitted bits. To corroborate the results, we compare the two proposed algorithms with the classical approach adopted in literature, defined by eq. \eqref{eq:ee_prob} \cite{Meshkati05, Buzzi08}.

\begin{figure*}[!htb]
    \begin{minipage}[t]{3 in}
        \begin{algorithm}[H]
        \caption{EE-SE with Verhulst Optimum Power Allocation}
        \begin{algorithmic}
            \REQUIRE $i \leftarrow 1$,\, It, \,\, $p_k[0]=\sigma^2_k, \,\forall k$
            \WHILE {$i\leq$ It}
                \FOR {$k=1:K$}
                    \STATE Evaluate $h_k$, $\widetilde{I}_k$;
                    \STATE Find $\gamma_k^*$ solving (\ref{eq:ee_se_deriv});
                    \STATE Find $p_k^*$ using Verhulst Algorithm.
                \ENDFOR
                \STATE $i \gets i + 1$
            \ENDWHILE
                \STATE Compute $\gamma_k$ for each user;
                \STATE Compute $K_{\rm out}$, where $k \in K_{\rm out}$ if $\gamma_{k} < \gamma_{k}^{*}$
                \IF {$\{K_{\rm out}\} \neq \emptyset$}
                    \STATE choose the user with worst channel gain in $K_{\rm out}$ ($j$-th user)
                	\STATE set $\gamma_j^*=0$;
                    \STATE go to the beginning.
                \ELSE
                    \RETURN $p^*_k \,\, \forall k$
                \ENDIF
        \end{algorithmic}
        \label{pcode:EE-SE-Verhulst}
        \end{algorithm}
    \end{minipage}
    \hfill
    \begin{minipage}[t]{3 in}
        \begin{algorithm}[H]
        \caption{EE-SE--$R_{k,\rm min}$ and Verhulst Power Allocation}
        \begin{algorithmic}
            \REQUIRE $i \leftarrow 1$,\, It, \,\, $p_k[0]=\sigma^2_k, \,\forall k$
            \STATE Compute $p_{k}^*$ as described in Algorithm \ref{pcode:EE-SE-Verhulst};
            \STATE Compute $\gamma_k$ for each user;
            \STATE Compute $K_{\rm out}$, where $k \in K_{\rm out}$ if $\gamma_{k} < \gamma_{k}^{*}$ and $r_k < R_{k,min}$
            \IF {$\{K_{\rm out}\} \neq \emptyset$}
                \STATE choose the user with worst channel gain in $K_{\rm out}$ ($j$-th user)
            	\STATE set $\gamma_j^*=0$;
                \STATE go to the beginning.
            \ELSE
                \RETURN $p^*_k \,\, \forall k$
            \ENDIF
        \end{algorithmic}
        \label{pcode:EE-SE-Verhulst_rate}
        \end{algorithm}
    \end{minipage}
    \hfill
\end{figure*}

\begin{figure}[!htbp]
    \centering
    \vspace{-2mm}
    a) \includegraphics[width=0.455\textwidth]{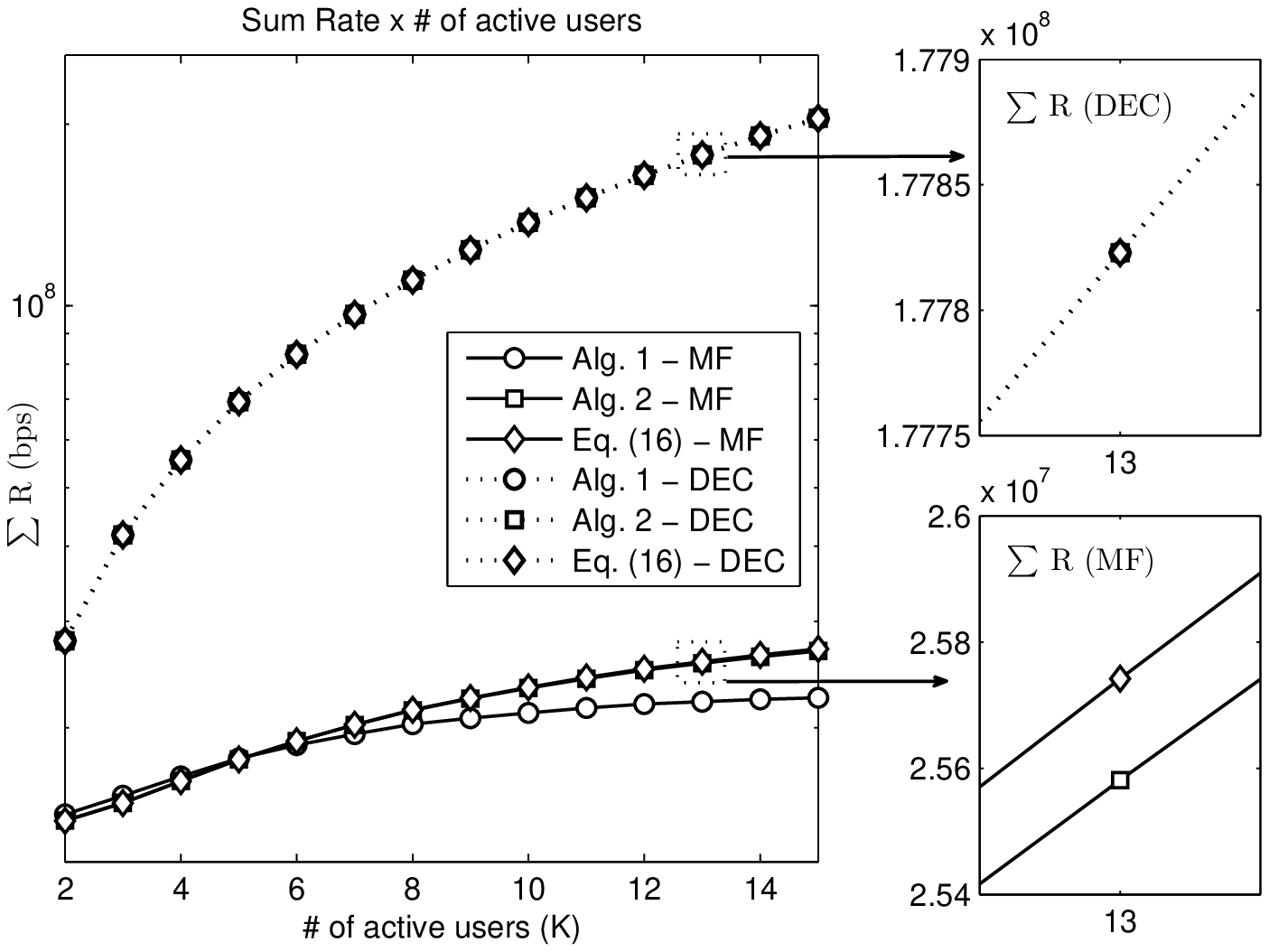}
    b) \includegraphics[width=0.455\textwidth]{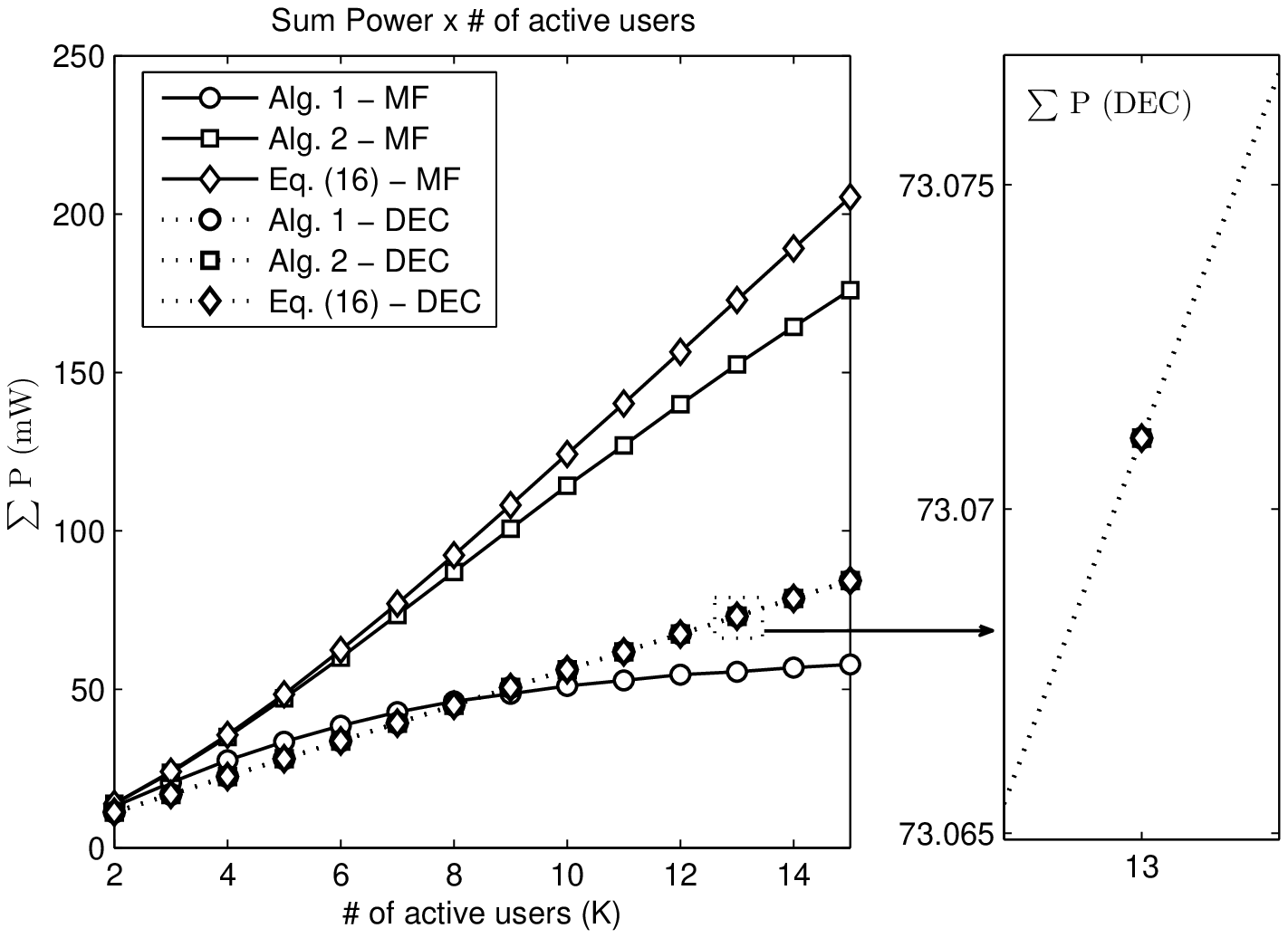}
    \caption{a) Sum rate (SR) and b) sum power (SP) for two different filters.}
    \vspace{-3mm}
\label{fig:Sum_RP}
\end{figure}

\begin{figure}[!htbp]
    \centering
    a) \includegraphics[width=0.455\textwidth]{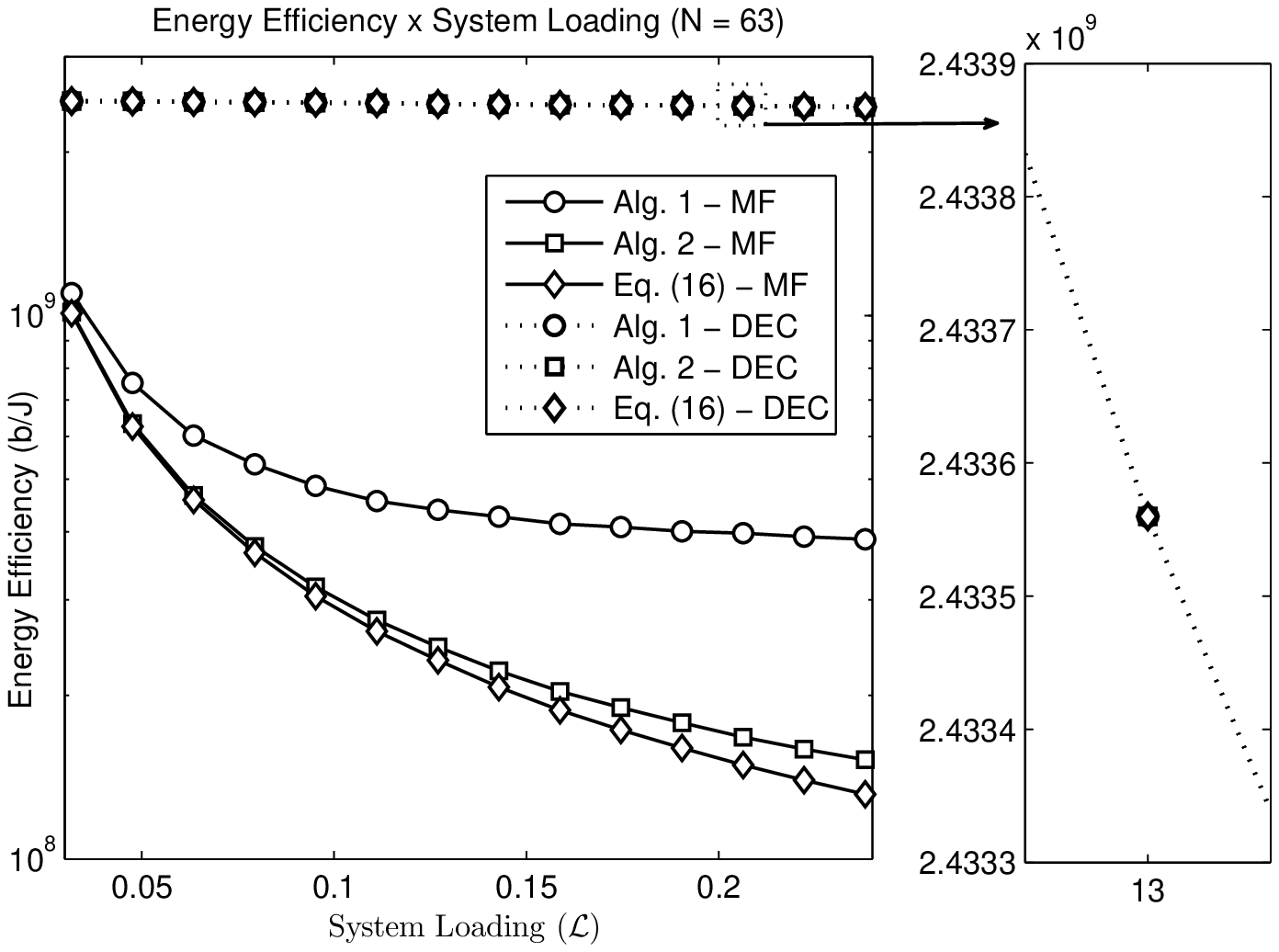}
    b) \includegraphics[width=0.455\textwidth]{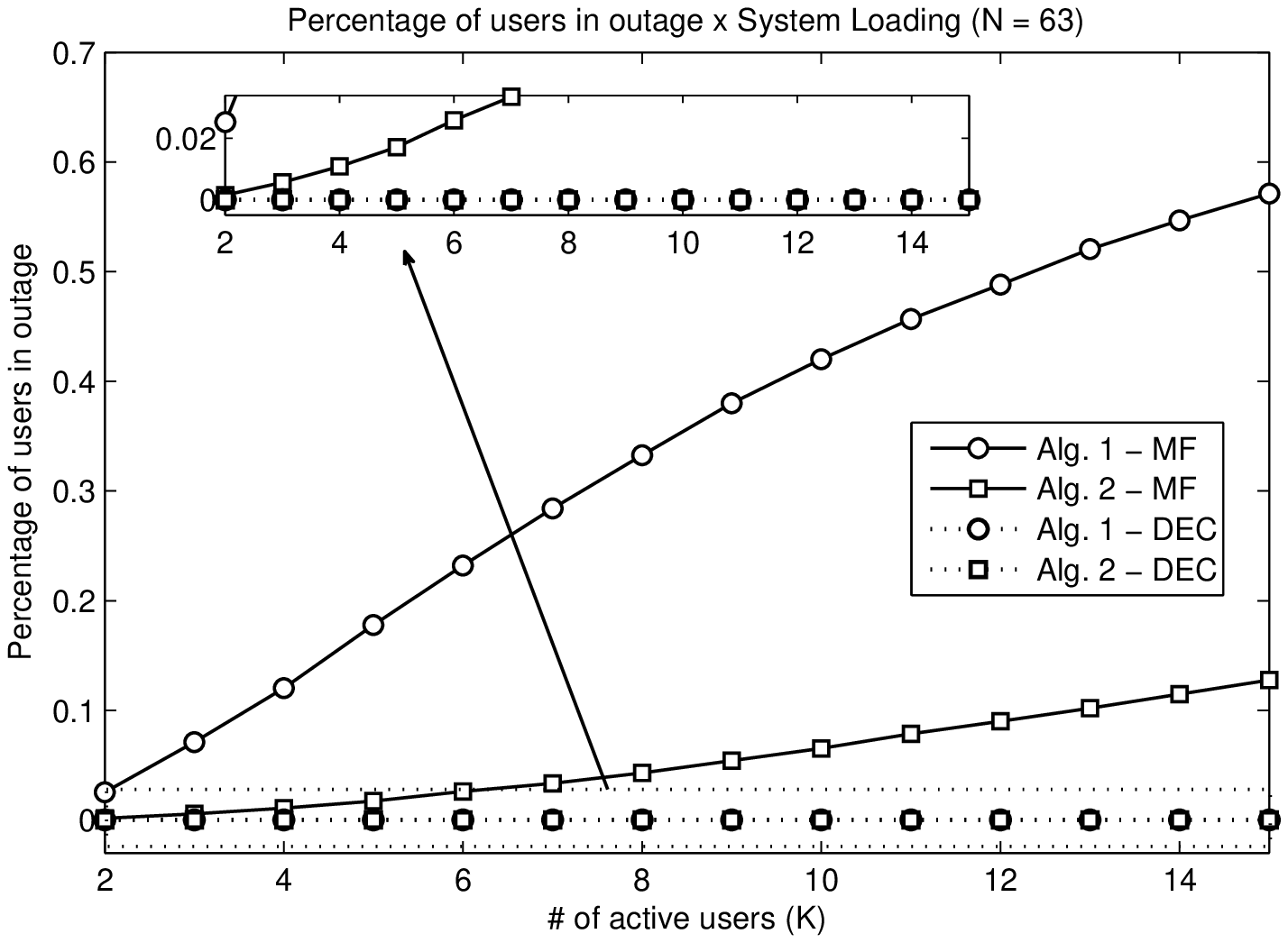}
    \caption{a) Energy Efficiency and b) outage probability for the two proposed algorithms and detectors.}
    \vspace{-3mm}
    \label{fig:EE_outage}
\end{figure}

Figs. \ref{fig:Sum_RP} and \ref{fig:EE_outage} bring four metric figures in order to analyze and to quantify performance gain of the two proposed algorithms, i.e., attainable sum rates of all users, $\sum R$, sum of power level consumption, including the circuit power, $\sum P$, the general energy efficiency (EE), obtained from the two algorithms for the two considered filters and, finally, the number of users put in outage.

From Fig. \ref{fig:Sum_RP} a) one can conclude that the problem defined in literature achieves the best result in terms of sum rate maximization, mainly when the system loading increases, since no user was put in outage, followed by Algorithm \ref{pcode:EE-SE-Verhulst_rate}, since a user is dropped only when the rate attainable by this user remain lower than the minimum rate value $R_{\min}$. As a consequence, the power consumption is increased, because users that don't achieve the optimum SINR try to achieve it using maximum power allowed, increasing remarkably the interference level. Since the decorrelator is more efficient than the matched filter on the MAI mitigation, the system is able to support more users under decorrelator multiuser filter. Besides, any sum power or sum rate performance difference between the two proposed algorithms could be noted, despite the evident efficiency increasing of both algorithms, since the achieved rate is larger and simultaneously the allocated power is smaller than those attainable with the matched filter.

Fig. \ref{fig:EE_outage} a) indicates the EE behavior against increasing system loading (when the number of active users increases). When $K\geq 7$, despite of the sum rate improvement obtained by the Algorithm \ref{pcode:EE-SE-Verhulst_rate} and the problem defined in eq. \eqref{eq:ee_prob} with matched filter over Algorithm \ref{pcode:EE-SE-Verhulst} with matched filter too, as one can see from Fig. \ref{fig:Sum_RP}, this improvement is obtained at cost of the system's energy efficiency degradation, Fig. \ref{fig:EE_outage}. This behavior is justified by the fact that there are users transmitting  with non-optimal powers in Algorithm \ref{pcode:EE-SE-Verhulst_rate} and in the problem defined by the literature. Besides, the best response in terms of EE is achieved by Algorithm \ref{pcode:EE-SE-Verhulst}, but incurs in more users in outage. Although there is a marginal power-rate trade-off difference among the two proposed algorithms, both are more efficient than the classical approach, mainly when system loading increases. As pointed out before, the multiuser decorrelating detector is more efficient than the MF for the two proposed algorithms, thanks to its improved capacity to provide multiple access interference mitigation.

Fig. \ref{fig:EE_outage} b) shows the impact of the MAI into receivers equipped with matched and decorrelator filters based systems. System loading was confined in the interval $\frac{K}{N}\in [ 0.0317;\,   0.2381]$. Hence, even under low system loading (around $20\%$ -- $24$\%), Algorithm \ref{pcode:EE-SE-Verhulst} based system is not able to achieve the maximum energy efficiency point for all users, since the required power to achieve the optimum SINR increases as the interference increases and then the maximum power available is overcame very soon. On the other hand, since Algorithm \ref{pcode:EE-SE-Verhulst_rate} allows users to transmit over a non-optimum power level scenario (as long as the minimum rate criterium is reached), the outage probability will be smaller. Again, thanks to the MAI mitigation characteristics of the decorrelator filter, the outage probability will be smaller than for MF. However, under extremely low loading (under 15\%), Algorithm \ref{pcode:EE-SE-Verhulst_rate} with MF presents lower outage probability than Algorithm \ref{pcode:EE-SE-Verhulst} with DEC filter; of course, with increasing system loading the outage probability becomes higher than Algorithm \ref{pcode:EE-SE-Verhulst} with DEC.

It is worth to note that the performance gaps among the two proposed algorithms -- mainly deployed with matched filter based systems -- can be explained by the numerical value for the minimum rate adopted, which requires a low spectral efficiency ($\eta_k = 0.5$) to be achieved, while allow a better visualization of the performance difference. Adopting a higher value for the minimum rate criterion, the expectation is that the outage probability will be increased for the two filters (MF and DEC), while the performance difference among the Algorithm \ref{pcode:EE-SE-Verhulst} and \ref{pcode:EE-SE-Verhulst_rate} will be decreased.

\begin{figure}[!htbp]
    \centering
    \vspace{-2mm}
    a) \includegraphics[width=0.455\textwidth]{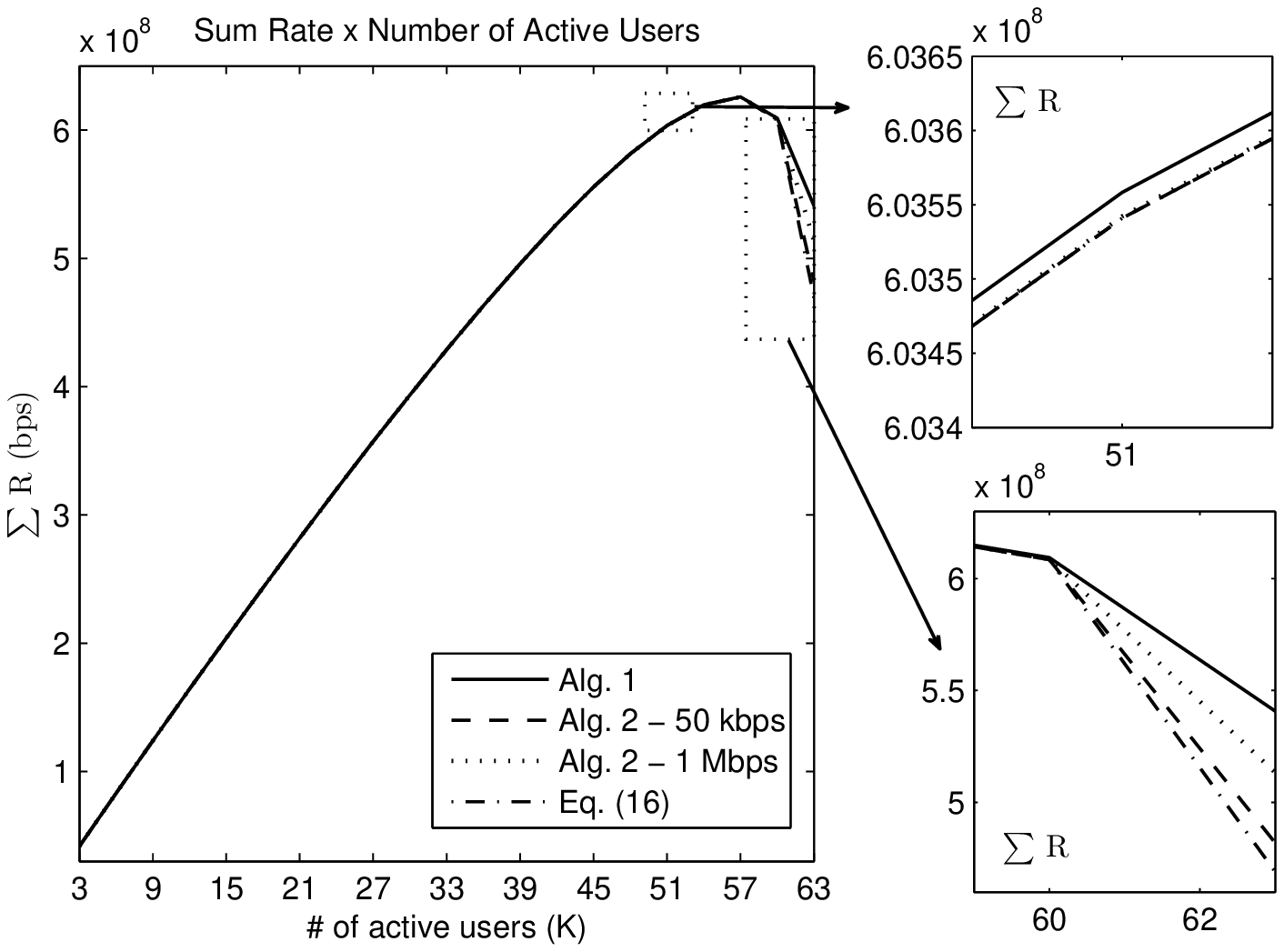}
    b) \includegraphics[width=0.455\textwidth]{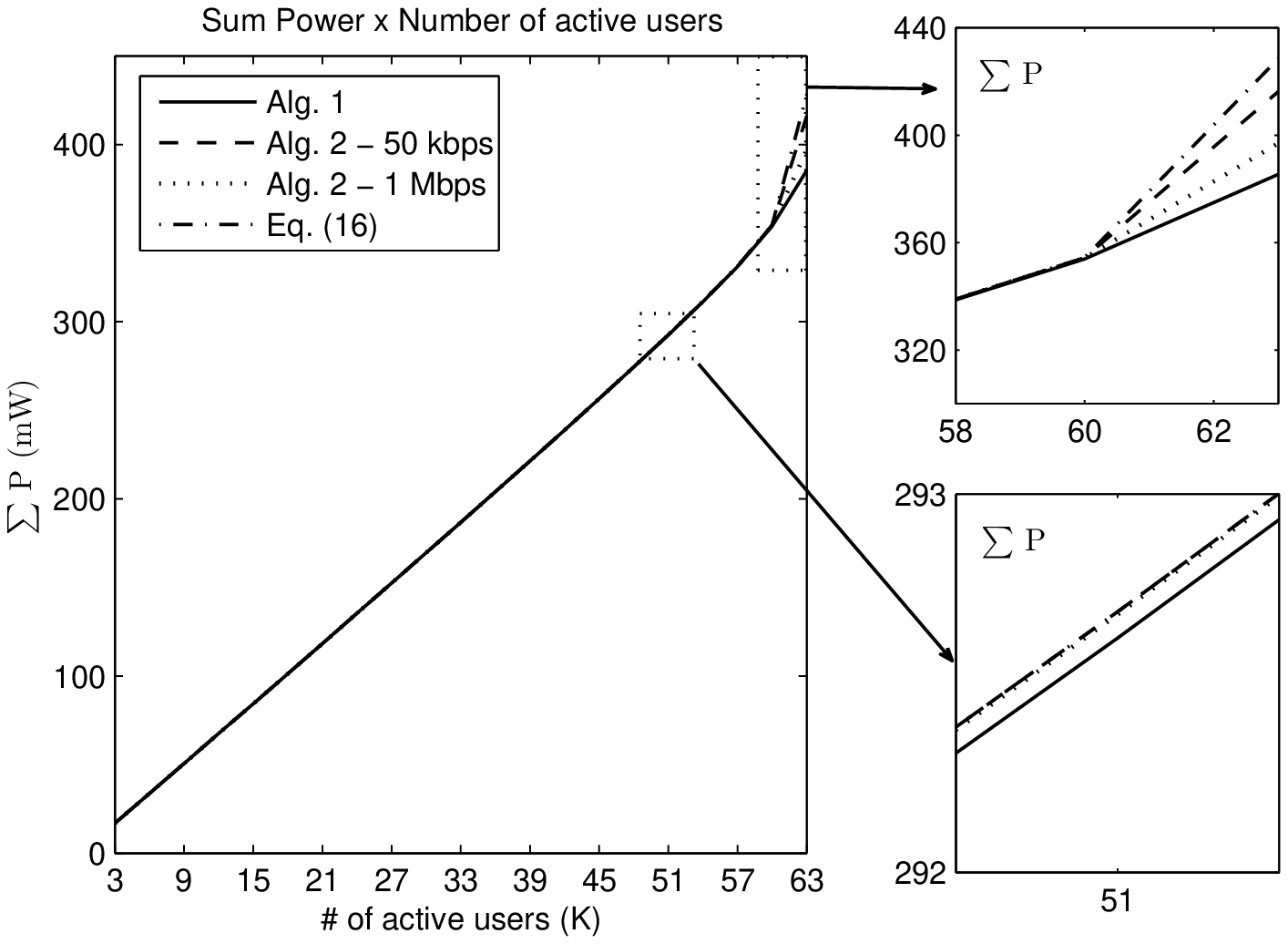}
    \caption{a) Sum rate (SR) and b) sum power (SP) for decorrelator filter.}
    \vspace{-3mm}
\label{fig:Sum_RP_dec}
\end{figure}

\begin{figure}[!htbp]
    \vspace{-2mm}
    \centering
    a) \includegraphics[width=0.465\textwidth]{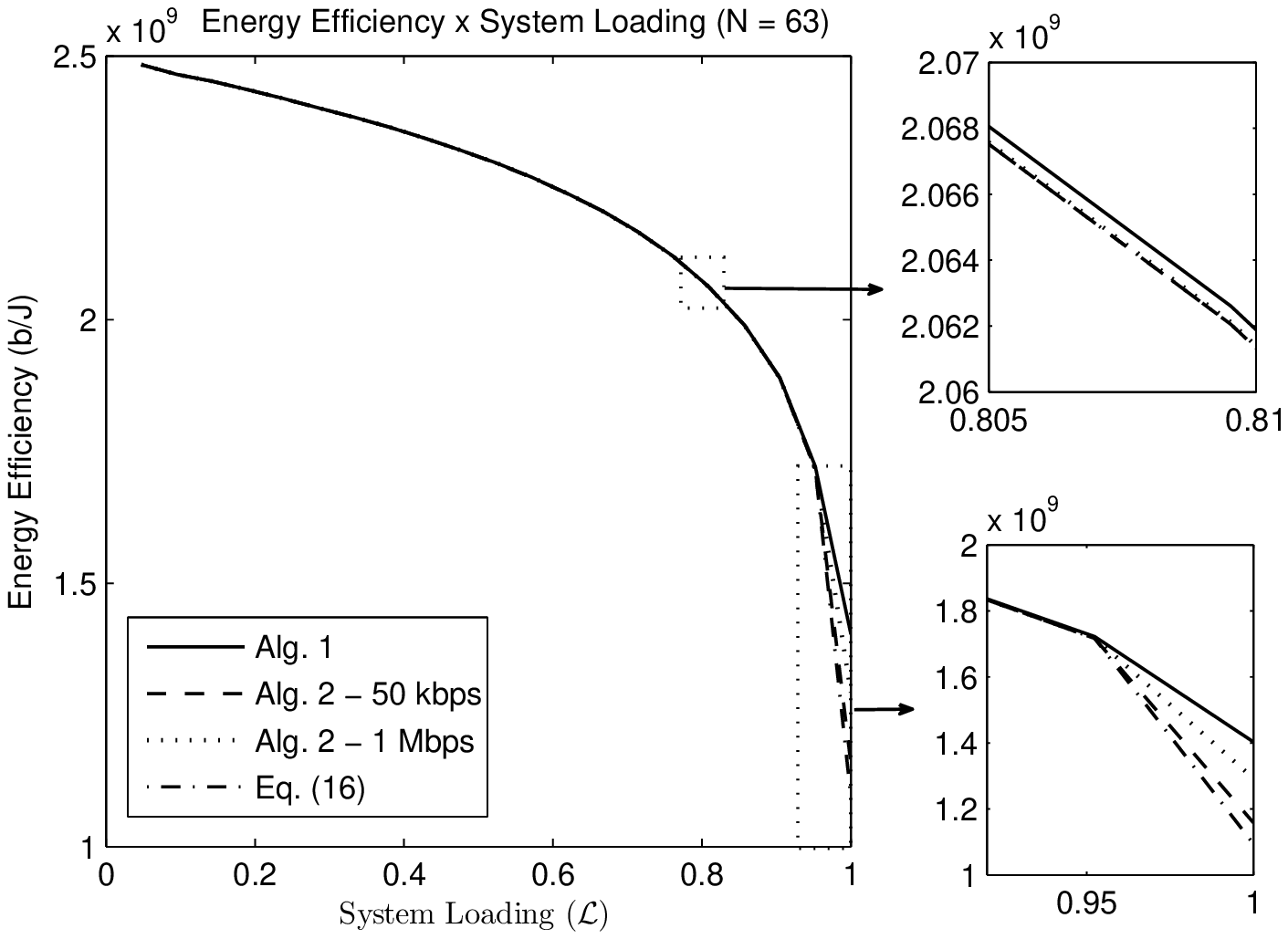}
    b) \includegraphics[width=0.465\textwidth]{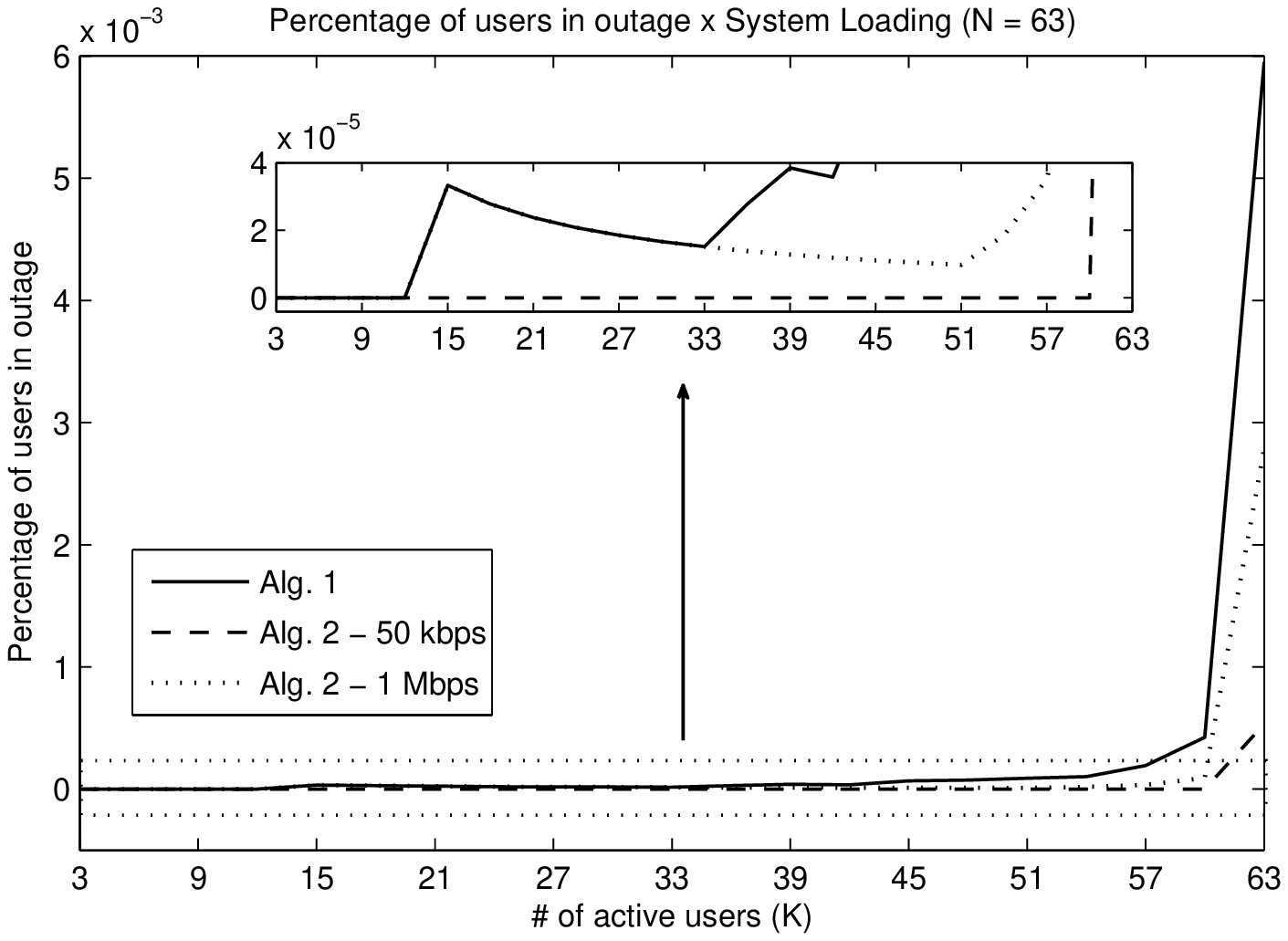}
    \caption{a) Energy Efficiency; b) outage probability for the two proposed algorithms with decorrelator.}
    \vspace{-3mm}
    \label{fig:EE_outage_dec}
\end{figure}

To corroborate the efficiency of decorrelator and the conclusions about rate criterium, we simulate the same metric figures (Figs. \ref{fig:Sum_RP_dec} and \ref{fig:EE_outage_dec}) now at full-loading ($K \in [3; 63]$) and with two different rate criteria ($R_{k, \min} = 50$ kbps and $R_{k, \min} = 1$ Mbps). For the first two metrics (sum-rate and sum-power, Fig. \ref{fig:Sum_RP_dec}) the results demonstrate that Algorithm \ref{pcode:EE-SE-Verhulst} obtains the best results in terms of sum-rate and sum-power, followed by Algorithm \ref{pcode:EE-SE-Verhulst_rate} and the algorithm described by eq. \eqref{eq:ee_prob}. Note that the bigger the minimum rate criterium, the closer Algorithm \ref{pcode:EE-SE-Verhulst_rate} results are from Algorithm \ref{pcode:EE-SE-Verhulst}.

From Fig. \ref{fig:EE_outage_dec} a), we can see that the achieved energy efficiency decreases when system loading increases, but still bigger than the obtained for MF in low system loading. Again, Algorithm \ref{pcode:EE-SE-Verhulst} is the most efficient, followed by Algorithm \ref{pcode:EE-SE-Verhulst_rate} and the literature common approach. For outage probability, Fig. \ref{fig:EE_outage_dec} b) demonstrate again that Algorithm \ref{pcode:EE-SE-Verhulst} presents the higher probability, but even at full loading this probability are lower than the obtained by MF at lower loadings. Again, higher the minimum rate, closer Algorithm \ref{pcode:EE-SE-Verhulst_rate} is from Algorithm \ref{pcode:EE-SE-Verhulst} and lower the minimum rate, closer Algorithm \ref{pcode:EE-SE-Verhulst_rate} is from literature's common approach.

\section{Conclusions} \label{sec:conclusions}
In this work we have analyzed the distributed energy efficiency (EE) cost function from the perspective of two conflicting metrics, throughput maximization and power level consumption minimization, as well as the impact of multiuser filter deployment over EE-SE trade-off.

We have found that SINR under the max-EE point equilibrium is almost the same whatever the level of multiple access interference becomes, mainly if interference level is medium or high. For MF, the best energy-spectral efficiencies trade-off consists in allocating to each node the necessary transmit power to achieve the best SINR response, which guarantees the maximal EE, while SE can be determined by the attainable rate in each node gives by the Shannon capacity equation.

Employing different figures of merit, numerical results indicated that deploying both proposed power allocation algorithms the linear multiuser filter is much more efficient than conventional matched filter receiver.

Finally, since the decorrelator detector is more efficient in providing MAI mitigation, new formulation for the max-EE \emph{versus} opt-SE trade-off problem, considering multi-objetive techniques would be proposed as a new research direction in the field.


\end{document}